\documentclass[a4j,10pt]{article}
\usepackage{amsmath}
\usepackage{amssymb}
\usepackage{amsthm}
\usepackage{graphics}
\usepackage[dvips]{graphicx}
\usepackage{amscd}

\numberwithin{equation}{section}
\theoremstyle{definition}
\newtheorem{example}{example}[section]
\newtheorem{df}[example]{Definition}
\theoremstyle{remark}

\theoremstyle{plain}
\newtheorem{lm}[example]{Lemma}

\newtheorem{theorem}{Theorem}[section]

\newtheorem{lemma}[theorem]{Lemma}

\newtheorem{remark}[theorem]{Remark}

\theoremstyle{remark}

\begin{document}

\title{Algebraic Ricci solitons of four-dimensional pseudo-Riemannian
generalized symmetric spaces}
\author{ W. Batat and K. Onda\thanks{%
First author was supported by funds of JSPS and ENSET d'Oran. The second
author was supported by funds of Nagoya University and OCAMI. \newline
} }
\maketitle

\begin{abstract}
We completely classify the algebraic Ricci solitons of four-dimensional
pseudo-Riemannian generalized symmetric spaces.
\end{abstract}

Wafaa Batat

\textsc{\'{E}}cole Normale Sup\'{e}rieure d'Enseignement Technologique d'Oran

B.P 1523 El M'naouar Oran 31000, Algeria

\textit{E-mail address}: wafa.batat@enset-oran.dz

\bigskip


Kensuke Onda

Graduate School of Mathematics, Nagoya University, Furocho, Chikusaku, \ \ 

Nagoya, Japan / POSTAL CODE: 464-860 

\textit{E-mail address}: kensuke.onda@math.nagoya-u.ac.jp

\textit{PHONE}: +81-52-789-2429

\textit{FACSIMILE}: +81-52-789-2829 


\section{\protect\LARGE Introduction and preliminaries }

\setcounter{equation}{0}

The concept of the Algebraic Ricci soliton was first introduced by Lauret in
Riemannian case (\cite{L01}). The definition extends to the
pseudo-Riemannian case:

\begin{df}
\textit{Let } $\left( \mathit{G,g}\right) $\textit{\ be a simply connected
Lie group equipped with the left-invariant pseudo-Riemannian metric $g$, and
let } ${\mathfrak{g}}$ denote \textit{the Lie algebra }of $\mathit{G.}$ 
\textit{Then $g$ is called an algebraic Ricci soliton if it satisfies } 
\begin{equation}
\mathrm{Ric}=c\mathrm{Id}+D  \label{Soliton}
\end{equation}%
\textit{where }$\mathrm{Ric}$\textit{\ denotes the Ricci operator, $c$ is a
real number, and }$\mathrm{D\in Der}\left( \mathfrak{g}\right) $ ($\mathrm{D}
$ is \textit{a derivation of }$\mathit{\mathfrak{g}}$\textit{)}, that is:%
\begin{equation}
\mathit{\mathrm{D}[X,Y]=[\mathrm{D}X,Y]+[X,\mathrm{D}Y]\ \,}\text{\ \textit{%
for any }}\mathit{\,\ X,Y\in \mathfrak{g}.}  \label{Der}
\end{equation}%
\textit{\ In particular, an algebraic Ricci soliton on a solvable Lie group,
(a nilpotent Lie group) is called a sol-soliton (a nil-soliton).}
\end{df}

Obviously, Einstein metrics on a Lie group are algebraic Ricci solitons.

A Ricci soliton metric $g$ on a manifold $M$ is a pseudo-Riemannian metric
satisfying 
\begin{equation}
\varrho \lbrack g]=cg+L_{X}g\ ,  \label{Ricci soliton}
\end{equation}%
where $L_{X}$ denotes the Lie derivative in the direction of the vector
field $X$, $\varrho $ is the Ricci tensor and $c$ is a real constant. A
Ricci soliton is said to be shrinking, steady or expanding, if $\lambda >0$, 
$\lambda =0$ or $\lambda <0$, respectively. Moreover, we say that a Ricci
soliton $(M,g)$ is a gradient Ricci soliton if it admits a vector field $X$
satisfying $X=\mathrm{grad}\left( h\right) $, for some potential function $h$
(see also \cite{CK04}).

Furthermore, the condition 
\eqref{Ricci soliton}
is equivalent to $g\left( t\right) =\left( -2ct+1\right) \varphi _{s\left(
t\right) }^{\ast }$ being a solution to the Ricci flow 
\begin{equation*}
\frac{\partial }{\partial t}g(t)_{ij}=-2\varrho \lbrack g(t)]_{ij}\ ,
\end{equation*}%
where $\varphi _{s}$ is the family of diffeomorphisms generated by $X$ with $%
s\left( t\right) =\frac{1}{c}\ln (-2ct+1).$

In \cite{L01}, Lauret studied the relation between sol-solitons and Ricci
solitons on Riemannian manifolds. More precisely, he proved that any
left-invariant Riemannian sol-soliton metric is a Ricci soliton. This was
extended by the second author to the pseudo-Riemannian case :

\begin{theorem}[\protect\cite{O11}]
\label{RS} Let $\left( G,g\right) $ be a simply connected Lie group endowed
with a left-invariant pseudo-Riemannian metric $g$. If $g$ is sol-soliton,
then $g$ is the Ricci soliton, that is, $g$ satisfy 
\eqref{Ricci soliton}%
, such that%
\begin{equation*}
X=\dfrac{d\varphi _{t}}{dt}\Big|_{t=0}(p)\text{ \ and }\mathrm{\exp }\left( 
\frac{t}{2}\mathrm{D}\right) =d\varphi _{t}|_{e},
\end{equation*}%
where $e$ denotes the identity element of $G$ .
\end{theorem}

Note that the above theorem is correct, if changing sol-soliton into an
algebraic Ricci soliton.

On the other hand, if $\left( M,g\right) $ is a homogeneous
(pseudo-)Riemannian manifold, then there exists a group $G$ of isometries
acting transitively on it \cite{O}. Such $(M,g)$ can be then identified with 
$(G/H,g)$, where $H$ is the isotropy group at a fixed point $p$ of $M$. Let $%
\mathfrak{g}$ denote the Lie algebra of $G$ and fix an \textrm{Ad}$\left(
H\right) $-invariant decomposition $\mathfrak{g}=\mathfrak{h}\oplus m$,
where $\mathfrak{h}$ is the Lie algebra of $H$. The space $m$ is naturally
identified with $T_{p}M.$ It must be noted, that such decomposition exists
always in the Riemannian case since homogeneous Riemannian manifolds are
reductive. However, in the pseudo-Riemannian setting the homogeneous
pseudo-Riemannian manifold need to be reductive homogeneous. In fact, a
three-dimensional homogeneous Lorentzian manifold is necessarily reductive.
This was proved in \cite{FR} and also follows independently from the
classification obtained by Calvaruso in \cite{C07}. Furthermore, the
existence of a non-reductive four-dimensional pseudo-Riemannian homogeneous
manifolds was proven in \cite{FR}. 

Homogeneous Ricci solitons have been investigated in \cite{Jac}. A natural
generalization of algebraic Ricci solitons on Lie groups to homogeneous
spaces is the following \cite{Jac}:

\begin{df}
Let $\left( M=G/H,g\right) $ be a homogeneous Riemannian manifold. \textit{%
Then $g$ is called an algebraic Ricci soliton if } 
\begin{equation}
\mathrm{Ric}=c\mathrm{Id}+pr\circ D  \label{AS}
\end{equation}%
\textit{where }$\mathrm{Ric}$\textit{\ denotes the Ricci operator of }$m$, $%
pr:\mathfrak{g}\rightarrow m$\textit{, $c$ is a real number, and }$\mathrm{%
D\in Der}\left( \mathfrak{g}\right) .$
\end{df}

Note that, the above definition can be extended to the pseudo-Riemannian
case, if changing homogeneous Riemannian manifold into reductive homogeneous
pseudo-Riemannian manifold. In \cite{BO}, we obtained the classification of
algebraic Ricci solitons of three-dimensional Lorentzian Lie groups.

In \cite{CK}, four-dimensional generalized symmetric spaces have been
completely classified. They are divided into four classes, named A, B, C and
D and the (pseudo-)Riemannian metrics can have any signature. All these
spaces are reductive homogeneous. In \cite{CD}, the Levi Civita connection,
the curvature tensor and the Ricci tensor of these spaces are computed. We
will use the results of these computations to study the algebraic Ricci
solitons of these spaces.

\section{\protect\LARGE On generalized symmetric spaces}

\setcounter{equation}{0}

We start by recalling the definition of generalized symmetric space. Let $%
\left( M,g\right) $ be a (pseudo-)Riemannian manifold. A \textit{regular }$s$%
\textit{-structure} on $M$ is a family of isometries $\left\{ s_{p}\text{ }|%
\text{ }p\in M\right\} $ of $\left( M,g\right) $ such that

\begin{itemize}
\item the mapping $M\times M\rightarrow M:\left( p,q\right) \mapsto
s_{p}\left( q\right) $ is smooth,

\item $\forall p\in M:p$ is an isolated fixed point of $s_{p}$,

\item $\forall p,q\in M:s_{p}\circ s_{q}=s_{s_{p}\left( q\right) }\circ
s_{p}.$
\end{itemize}

$s_{p}$ is called a symmetry centered at $p.$ The \textit{order }of a regular%
\textit{\ }$s$-structure is the smallest integer $m\geqslant 2$ such that $%
s_{p}^{m}=\mathrm{id}_{M}$ for all $p\in M.$ If such an integer does not
exist, we say that the regular $s$-structure has order infinity. A \textit{%
generalized symmetric space} is a connected, pseudo-Riemannian manifold,
carrying at least one regular $s$-structure. In particular, a generalized
symmetric space is a symmetric space if and only if it admits a regular $s$%
-structure of order 2. The \textit{order} of a generalized symmetric space
is the minimum of orders of all possible $s$-structures on it. Furthermore,
if $\left( M,g\right) $ is a generalized symmetric space then it is
homogeneous, that is, the full isometry group $I\left( M\right) $ of $M$
acts transitively on it, this means that $\left( M,g\right) $ can be
identified with $\left( G/H,g\right) $, where $G\subset I\left( M\right) $
is a subgroup of $I\left( M\right) $ acting transitively on $M$ and $\ H$ is
the isotropy group at a fixed point $o\in M.$

Generalized symmetric spaces of low dimension have been completely
classified. The following Theorem presents a complete classification of
four-dimensional pseudo-Riemannian generalized symmetric spaces.

\begin{theorem}[\protect\cite{CK}]
All non-symmetric, simply connected generalized symmetric spaces $\left(
M,g\right) $\ of dimension 4 are of order 3 or 4, or infinity. All these
spaces are indecomposable, and belong, up to isometry, to the following four
types:

\begin{itemize}
\item $\mathrm{Type}$ $\mathrm{A.}$ The underlying homogeneous space is $G/H$%
, where%
\begin{equation*}
G=\left( 
\begin{array}{ccc}
a & b & u \\ 
c & d & v \\ 
0 & 0 & 1%
\end{array}%
\right) ,\text{ \ }H=\left( 
\begin{array}{ccc}
\cos t & -\sin t & 0 \\ 
\sin t & \cos t & 0 \\ 
0 & 0 & 1%
\end{array}%
\right)
\end{equation*}%
with $ad-bc=1.$ $\left( M,g\right) $ is the space $%
\mathbb{R}
^{4}\left( x,y,u,v\right) $ with the pseudo-Riemannian metric%
\begin{eqnarray}
g &=&\lambda \left[ \left( 1+y^{2}\right) dx^{2}+\left( 1+x^{2}\right)
dy^{2}-2x y dx dy\right] /\left( 1+x^{2}+y^{2}\right)  \label{gApR} \\
&&\pm \left[ \left( -x+\sqrt{1+x^{2}+y^{2}}\right) du^{2}+\left( x+\sqrt{%
1+x^{2}+y^{2}}\right) dv^{2}-2y^{2}dudv\right] ,  \notag
\end{eqnarray}%
where $\lambda \neq 0$ is a real constant. The order is $k=3$ and possible
signature are $\left( 4,0\right) ,\left( 0,4\right) ,\left( 2,2\right) .$

\item $\mathrm{Type}$ $\mathrm{B.}$ The underlying homogeneous space is $G/H$%
, where%
\begin{equation*}
G=\left( 
\begin{array}{cccc}
e^{-\left( x+y\right) } & 0 & 0 & a \\ 
0 & e^{x} & 0 & b \\ 
0 & 0 & e^{y} & c \\ 
0 & 0 & 0 & 1%
\end{array}%
\right) ,\text{ \ }H=\left( 
\begin{array}{cccc}
1 & 0 & 0 & -w \\ 
0 & 1 & 0 & -2w \\ 
0 & 0 & 1 & 2w \\ 
0 & 0 & 0 & 1%
\end{array}%
\right) .
\end{equation*}%
$\left( M,g\right) $ is the space $%
\mathbb{R}
^{4}\left( x,y,u,v\right) $ with the pseudo-Riemannian metric%
\begin{equation}
g=\lambda \left( dx^{2}+dy^{2}+dxdy\right) +e^{-y}\left( 2dx+dy\right)
dv+e^{-x}\left( dx+2dy\right) du  \label{gB}
\end{equation}%
where $\lambda $ is a real constant. The order is $k=3$ and the signature is
always $\left( 2,2\right) .$

\item $\mathrm{Type}$ $\mathrm{C.}$ The underlying homogeneous space $G/H$
is the matrix group%
\begin{equation*}
G=\left( 
\begin{array}{cccc}
e^{-t} & 0 & 0 & x \\ 
0 & e^{t} & 0 & y \\ 
0 & 0 & 1 & z \\ 
0 & 0 & 0 & 1%
\end{array}%
\right) .
\end{equation*}%
$\left( M,g\right) $ is the space $%
\mathbb{R}
^{4}\left( x,y,z,t\right) $ with the pseudo-Riemannian metric%
\begin{equation}
g=\pm \left( e^{2t}dx^{2}+e^{-2t}dy^{2}\right) +dzdt.  \label{gC}
\end{equation}%
The order is $k=3$ and possible signature are $\left( 1,3\right) ,\left(
3,1\right) .$

\item $\mathrm{Type}$ $\mathrm{D.}$ The underlying homogeneous space is $G/H$
where%
\begin{equation*}
G=\left( 
\begin{array}{ccc}
a & b & x \\ 
c & d & y \\ 
0 & 0 & 1%
\end{array}%
\right) ,\text{ \ }H=\left( 
\begin{array}{ccc}
e^{t} & 0 & 0 \\ 
0 & e^{-t} & 0 \\ 
0 & 0 & 1%
\end{array}%
\right)
\end{equation*}%
with $ad-bc=1.$ $\left( M,g\right) $ is the space $%
\mathbb{R}
^{4}\left( x,y,u,v\right) $ with the pseudo-Riemannian metric%
\begin{eqnarray}
g &=&-2\cosh \left( 2u\right) \cos \left( 2v\right) dxdy+\lambda \left(
du^{2}-\cosh ^{2}\left( 2u\right) dv^{2}\right)  \label{gD} \\
&&+\left( \sinh \left( 2u\right) -\cosh \left( 2u\right) \sin \left(
2v\right) \right) dx^{2}+\left( \sinh \left( 2u\right) +\cosh \left(
2u\right) \sin \left( 2v\right) \right) dy^{2},  \notag
\end{eqnarray}%
where $\lambda \neq 0$ is a real constant. The order is infinite and the
signature is $\left( 2,2\right) .$
\end{itemize}
\end{theorem}

\section{\protect\LARGE Algebraic Ricci soliton of spaces of type A with
neutral signature}

\setcounter{equation}{0}

Let $\left( M,g\right) $ be a four-dimensional generalized pseudo-Riemannian
symmetric space and denote by $\nabla $ and $R$ the Levi-Civita connection
and\ the Riemann curvature tensor of $M$ respectively. Throughout this
paper, we will always use the sign convention%
\begin{equation*}
R\left( X,Y\right) =\nabla _{\left[ X,Y\right] }-\left[ \nabla _{X},\nabla
_{Y}\right] .
\end{equation*}%
The Ricci tensor of $\left( M,g\right) $ is defined by%
\begin{equation*}
\varrho \left( X,Y\right) =\overset{4}{\underset{k=1}{\sum }}\varepsilon
_{k}g(R(X,e_{k})Y,e_{k}),
\end{equation*}%
where $\{e_{1},e_{2},e_{3},e_{4}\}$ is a pseudo-orthonormal frame field,
with $g\left( e_{k},e_{k}\right) =\varepsilon _{k}=\pm 1.$ The Ricci
operator \textrm{Ric} is then given by%
\begin{equation*}
\varrho \left( X,Y\right) =g\left( \mathrm{Ric}\left( X\right) ,Y\right) .
\end{equation*}

Now, consider a four-dimensional generalized symmetric space $\left(
M=G/H,g\right) $ of type A, with signature $\left( 2,2\right) .$ Then,
taking into account the results of \cite{CK} and \cite{DM} the Lie algebra $%
\mathfrak{g}$ of the Lie group $G\ $may be decomposed into vector spaces
direct sum $\mathfrak{g}$ $=\mathfrak{h\oplus }m$ where $\mathfrak{h}$
denotes the Lie algebra of $H$ and $m$ is a vector space of $\mathfrak{g}.$

The Lie algebra $\mathfrak{g}$\ admits a basis $\left\{
U_{1},U_{2},U_{3},U_{4},U_{5}\right\} $, where $\left\{
U_{1},U_{2},U_{3},U_{4}\right\} $ is an orthogonal basis of $m$ and $\left\{
U_{5}\right\} $ basis of $\mathfrak{h,}$ such that the Lie bracket $\left[ 
\text{ },\right] $ on $\mathfrak{g}$ and the scalar product $\left\langle 
\text{ },\right\rangle $ on $m$ are given, respectively, by%
\begin{equation}
\begin{array}[t]{c|ccccc}
\left[ \text{ },\right] & U_{1} & U_{2} & U_{3} & U_{4} & U_{5} \\ \hline
U_{1} & 0 & 0 & -\delta U_{1} & \delta U_{2} & U_{2} \\ 
U_{2} & 0 & 0 & \delta U_{2} & \delta U_{1} & -U_{1} \\ 
U_{3} & \delta U_{1} & -\delta U_{2} & 0 & -2\delta ^{2}U_{5} & -2U_{4} \\ 
U_{4} & -\delta U_{2} & -\delta U_{1} & 2\delta ^{2}U_{5} & 0 & 2U_{3} \\ 
U_{5} & -U_{2} & U_{1} & 2U_{4} & -2U_{3} & 0%
\end{array}
\label{LieA}
\end{equation}%
where $\delta >0$ is a real constant, and 
\begin{equation*}
\begin{array}{c|cccc}
\left\langle \text{ },\right\rangle & U_{1} & U_{2} & U_{3} & U_{4} \\ \hline
U_{1} & 1 & 0 & 0 & 0 \\ 
U_{2} & 0 & 1 & 0 & 0 \\ 
U_{3} & 0 & 0 & -2 & 0 \\ 
U_{4} & 0 & 0 & 0 & -2%
\end{array}%
\end{equation*}

We start by recalling the following result on the curvature tensor and the
Ricci tensor of four-dimensional generalized symmetric spaces of type A (see 
\cite{CD}).

\begin{lm}
\label{ApR} Let $M$ be a four-dimensional generalized symmetric space of
type A, with signature $\left( 2,2\right) .$ Then, there exist a
pseudo-orthonormal frame field $\left\{ e_{1}=U_{1},e_{2}=U_{2},e_{3}=\frac{1%
}{\sqrt{2}}U_{3},e_{4}=\frac{1}{\sqrt{2}}U_{4}\right\} $ on $M$, with $%
\left\langle e_{1},e_{1}\right\rangle =\left\langle e_{2},e_{2}\right\rangle
=-\left\langle e_{3},e_{3}\right\rangle =-\left\langle
e_{4},e_{4}\right\rangle =1.$ The non-vanishing components of the
Levi-Civita connection $\nabla $ of $M$ are given by%
\begin{eqnarray*}
\nabla _{e_{1}}e_{1} &=&-\frac{\delta }{\sqrt{2}}e_{3},\text{ \ }\nabla
_{e_{1}}e_{2}=\frac{\delta }{\sqrt{2}}e_{4},\text{ \ }\nabla _{e_{1}}e_{3}=-%
\frac{\delta }{\sqrt{2}}e_{1},\text{ \ }\nabla _{e_{1}}e_{4}=\frac{\delta }{%
\sqrt{2}}e_{2},\text{\ } \\
\nabla _{e_{2}}e_{1} &=&\frac{\delta }{\sqrt{2}}e_{4},\text{ \ \ \ }\nabla
_{e_{2}}e_{2}=\frac{\delta }{\sqrt{2}}e_{3},\text{ \ }\nabla _{e_{2}}e_{3}=%
\frac{\delta }{\sqrt{2}}e_{2},\text{ \ \ \ }\nabla _{e_{2}}e_{4}=\frac{%
\delta }{\sqrt{2}}e_{1},
\end{eqnarray*}%
The only non-zero components of the Riemann curvature tensor $R\left(
X,Y,Z,W\right) =g\left( R\left( X,Y\right) Z,W\right) ,$ with respect to $%
\left\{ e_{1},e_{2},e_{3},e_{4}\right\} ,$ are 
\begin{eqnarray*}
R_{1212} &=&-R_{1234}=-\delta ^{2}, \\
R_{1313} &=&-R_{1324}=R_{1414}=R_{1423}=R_{2323}=R_{2424}=-\frac{\delta ^{2}%
}{2}.
\end{eqnarray*}%
Consequently, the non-zero components of the Ricci tensor are given by%
\begin{equation*}
\varrho \left( e_{3},e_{3}\right) =\varrho \left( e_{4},e_{4}\right)
=-\delta ^{2}.
\end{equation*}
\end{lm}

Now, let $\mathrm{D}\in \mathrm{Der}\left( \mathfrak{g}\right) $ where $%
\mathfrak{g}$ is the Lie algebra used in 
\eqref{LieA}%
. Put 
\begin{equation*}
\mathrm{D}U_{l}=\lambda _{l}^{1}U_{1}+\lambda _{l}^{2}U_{2}+\lambda
_{l}^{3}U_{3}+\lambda _{l}^{4}U_{4}+\lambda _{l}^{5}U_{5}\text{ \ for all }%
l=1,..,5.
\end{equation*}

Starting from 
\eqref{LieA}%
, we can write down 
\eqref{Der}
and we get%
\begin{eqnarray}
&&%
\begin{array}{c}
\lambda _{3}^{5}+\delta \left( 2\lambda _{1}^{2}+\lambda _{3}^{4}\right) =0,%
\end{array}
\label{sysDA} \\
&&%
\begin{array}{c}
\lambda _{3}^{5}+\delta \left( 2\lambda _{2}^{1}-\lambda _{3}^{4}\right) =0,%
\end{array}
\notag \\
&&%
\begin{array}{c}
\lambda _{4}^{5}+\delta \left( \lambda _{1}^{1}-\lambda _{2}^{2}\right) =0,%
\end{array}
\notag \\
&&%
\begin{array}{c}
\lambda _{2}^{1}-\lambda _{1}^{2}+\lambda _{4}^{3}=0,%
\end{array}
\notag \\
&&%
\begin{array}{c}
\lambda _{1}^{1}-\lambda _{2}^{2}+\delta \lambda _{5}^{4}=0,%
\end{array}
\notag \\
&&%
\begin{array}{c}
\lambda _{1}^{2}+\lambda _{2}^{1}+\delta \lambda _{5}^{3}=0,%
\end{array}
\notag \\
&&%
\begin{array}{c}
\lambda _{3}^{2}+\lambda _{4}^{1}+2\delta \lambda _{5}^{1}=0,%
\end{array}
\notag \\
&&%
\begin{array}{c}
\lambda _{3}^{1}-\lambda _{4}^{2}+2\delta \lambda _{5}^{2}=0,%
\end{array}
\notag \\
&&%
\begin{array}{c}
2\lambda _{4}^{1}-\lambda _{3}^{2}+\delta \lambda _{5}^{1}=0,%
\end{array}
\notag \\
&&%
\begin{array}{c}
2\lambda _{4}^{2}+\lambda _{3}^{1}-\delta \lambda _{5}^{2}=0,%
\end{array}
\notag \\
&&%
\begin{array}{c}
\lambda _{4}^{1}-2\lambda _{3}^{2}-\delta \lambda _{5}^{1}=0,%
\end{array}
\notag \\
&&%
\begin{array}{c}
\lambda _{4}^{2}+2\lambda _{3}^{1}+\delta \lambda _{5}^{2}=0,%
\end{array}
\notag \\
&&%
\begin{array}{c}
\lambda _{3}^{5}=\delta ^{2}\lambda _{5}^{3},\text{ \ }\lambda
_{4}^{3}=-\lambda _{3}^{4},\text{ \ }\lambda _{4}^{5}=\delta ^{2}\lambda
_{5}^{4},%
\end{array}
\notag \\
&&%
\begin{array}{c}
\lambda _{1}^{3}=\lambda _{1}^{4}=\lambda _{1}^{5}=\lambda _{2}^{3}=\lambda
_{2}^{4}=\lambda _{2}^{5}=\lambda _{3}^{3}=\lambda _{4}^{4}=\lambda
_{5}^{5}=0.%
\end{array}
\notag
\end{eqnarray}%
A standard computation proves that all solutions of 
\eqref{sysDA}
are given by%
\begin{eqnarray*}
&&%
\begin{array}{c}
\lambda _{2}^{1}=-\lambda _{1}^{2}-\delta \lambda _{5}^{3},%
\end{array}%
\begin{array}{c}
\lambda _{2}^{2}=\lambda _{1}^{1}+\delta \lambda _{5}^{4},%
\end{array}%
\begin{array}{c}
\lambda _{3}^{1}=-\lambda _{4}^{2}=-\delta \lambda _{5}^{2},%
\end{array}
\\
&&%
\begin{array}{c}
\lambda _{3}^{2}=\lambda _{4}^{1}=-\delta \lambda _{5}^{1},%
\end{array}%
\begin{array}{c}
\lambda _{4}^{3}=-\lambda _{3}^{4}=2\lambda _{1}^{2}+\delta \lambda _{5}^{3}.%
\end{array}%
\end{eqnarray*}%
So, we proved the following.

\begin{lemma}
Let $\mathfrak{g}=\mathfrak{h\oplus }m$ be the Lie algebra used in 
\eqref{LieA}%
. Then $\mathrm{D}\in \mathrm{Der}\left( \mathfrak{g}\right) $ if and only if%
\begin{equation*}
\mathrm{D}=\left( 
\begin{array}{ccccc}
\lambda _{1}^{1} & -\lambda _{1}^{2}-\delta \lambda _{5}^{3} & -\delta
\lambda _{5}^{2} & -\delta \lambda _{5}^{1} & \lambda _{5}^{1} \\ 
\lambda _{1}^{2} & \lambda _{1}^{1}+\delta \lambda _{5}^{4} & -\delta
\lambda _{5}^{1} & \delta \lambda _{5}^{2} & \lambda _{5}^{2} \\ 
0 & 0 & 0 & 2\lambda _{1}^{2}+\delta \lambda _{5}^{3} & \lambda _{5}^{3} \\ 
0 & 0 & -2\lambda _{1}^{2}-\delta \lambda _{5}^{3} & 0 & \lambda _{5}^{4} \\ 
0 & 0 & \delta ^{2}\lambda _{5}^{3} & \delta ^{2}\lambda _{5}^{4} & 0%
\end{array}%
\right)
\end{equation*}
\end{lemma}

Using the above lemma, we now prove the following.

\begin{theorem}
Let $\left( M=G/H,g\right) $ be a four-dimensional generalized symmetric
space of type A, with signature $\left( 2,2\right) $. Then, $M$ is an
algebraic Ricci solitons. In particular%
\begin{equation*}
pr\circ \mathrm{D=}\left( 
\begin{array}{cccc}
-\delta ^{2} & 0 & 0 & 0 \\ 
0 & -\delta ^{2} & 0 & 0 \\ 
0 & 0 & 0 & 0 \\ 
0 & 0 & 0 & 0%
\end{array}%
\right) \text{ and }c=\delta ^{2}.
\end{equation*}
\end{theorem}

\textbf{Proof. }Using Lemma \ref{ApR} we obtain that the Ricci operator of $%
\left( M=G/H,g\right) $ is given, with respect to the basis $\left\{
U_{1},U_{2},U_{3},U_{4},U_{5}\right\} ,$ by%
\begin{equation*}
\mathrm{Ric}=\left( 
\begin{array}{cccc}
0 & 0 & 0 & 0 \\ 
0 & 0 & 0 & 0 \\ 
0 & 0 & \delta ^{2} & 0 \\ 
0 & 0 & 0 & \delta ^{2}%
\end{array}%
\right) .
\end{equation*}%
Hence, the algebraic Ricci soliton condition 
\eqref{AS}
on $M$ is satisfied if and only if:%
\begin{eqnarray}
&&%
\begin{array}{c}
\lambda _{1}^{1}=-c=-\delta ^{2},%
\end{array}
\notag \\
&&%
\begin{array}{c}
\lambda _{1}^{2}=\lambda _{5}^{1}=\lambda _{5}^{2}=\lambda _{5}^{3}=\lambda
_{5}^{4}=0.\square%
\end{array}
\notag
\end{eqnarray}

\begin{remark}
Algebraic Ricci solitons of spaces of type A, the Riemannian case, is
obtained by the following. If we change $U_{3}$ and $U_{4}$ of
pseudo-Riemannian case to $\dfrac{1}{\delta } U_{3}$ and $\dfrac{1}{\delta }
U_{4}$ and put $\rho = -\delta ^2$, we obtain the Riemannian case of the Lie
bracket. So, it is easy to check that this case has a soliton.
\end{remark}

\section{\protect\LARGE Algebraic Ricci soliton of spaces of type B }

\setcounter{equation}{0}

Let $\left( M=G/H,g\right) $ be a four-dimensional generalized symmetric
space of type B, with signature $\left( 2,2\right) .$ Then, $\mathfrak{g}=%
\mathfrak{h\oplus }m$ and $\left\{ U_{1},U_{2},U_{3},U_{4}\right\} $ and $%
\left\{ U_{5}\right\} $ are respectively a basis of $m$ and $\mathfrak{h,}$
such that the Lie bracket $\left[ \text{ },\right] $ on $\mathfrak{g}$ and
the scalar product $\left\langle \text{ },\right\rangle $ on $m$ are given,
respectively, by%
\begin{equation}
\begin{array}{c|ccccc}
\left[ \text{ },\right] & U_{1} & U_{2} & U_{3} & U_{4} & U_{5} \\ \hline
U_{1} & 0 & 0 & -U_{1} & \varepsilon U_{5}+U_{2} & 0 \\ 
U_{2} & 0 & 0 & -\varepsilon U_{5}+U_{2} & U_{1} & 0 \\ 
U_{3} & U_{1} & \varepsilon U_{5}-U_{2} & 0 & 0 & 2U_{2} \\ 
U_{4} & -\varepsilon U_{5}-U_{2} & -U_{1} & 0 & 0 & -2U_{1} \\ 
U_{5} & 0 & 0 & -2U_{2} & 2U_{1} & 0%
\end{array}
\label{LieB}
\end{equation}%
where $\varepsilon =\pm 1$, and 
\begin{equation*}
\begin{array}{c|cccc}
\left\langle \text{ },\right\rangle & U_{1} & U_{2} & U_{3} & U_{4} \\ \hline
U_{1} & 0 & 0 & -1 & 0 \\ 
U_{2} & 0 & 0 & 0 & -1 \\ 
U_{3} & -1 & 0 & 2\lambda & 0 \\ 
U_{4} & 0 & -1 & 0 & 2\lambda%
\end{array}%
\end{equation*}

The following result was proven in \cite{CD}.

\begin{lm}
\label{B} Let $M$ be a four-dimensional generalized symmetric space of type
B, with signature $\left( 2,2\right) .$ Then, there exist a
pseudo-orthonormal frame field 
\begin{eqnarray*}
e_{1} &=&\left( \lambda -\frac{1}{2}\right) U_{1}+U_{2},\text{ \ \ \ \ \ }%
e_{2}=\left( \lambda -\frac{1}{2}\right) U_{3}+U_{4}, \\
e_{3} &=&\left( \lambda +\frac{1}{2}\right) U_{1}+U_{2},\text{ \ \ \ \ \ }%
e_{4}=\left( \lambda +\frac{1}{2}\right) U_{3}+U_{4},
\end{eqnarray*}%
on $M$, with $\left\langle e_{1},e_{1}\right\rangle =\left\langle
e_{2},e_{2}\right\rangle =-\left\langle e_{3},e_{3}\right\rangle
=-\left\langle e_{4},e_{4}\right\rangle =1.$ The Levi-Civita connection $%
\nabla $ of $M$ is determined by%
\begin{eqnarray*}
&&%
\begin{array}{c}
\nabla _{e_{1}}e_{1}=-e_{3},%
\end{array}%
\text{ \ \ }%
\begin{array}{c}
\nabla _{e_{2}}e_{1}=e_{4},%
\end{array}%
\text{ \ \ }%
\begin{array}{c}
\nabla _{e_{3}}e_{1}=-e_{3},%
\end{array}%
\text{ \ \ }%
\begin{array}{c}
\nabla _{e_{4}}e_{1}=e_{4},%
\end{array}
\\
&&%
\begin{array}{c}
\nabla _{e_{1}}e_{2}=e_{4},%
\end{array}%
\text{ \ \ \ \ }%
\begin{array}{c}
\nabla _{e_{2}}e_{2}=e_{3},%
\end{array}%
\text{ \ \ }%
\begin{array}{c}
\nabla _{e_{3}}e_{2}=e_{4},%
\end{array}%
\text{ \ \ \ \ \ }%
\begin{array}{c}
\nabla _{e_{4}}e_{2}=e_{3},%
\end{array}
\\
&&%
\begin{array}{c}
\nabla _{e_{1}}e_{3}=-e_{1},%
\end{array}%
\text{ \ \ }%
\begin{array}{c}
\nabla _{e_{2}}e_{3}=e_{2},%
\end{array}%
\text{ \ \ }%
\begin{array}{c}
\nabla _{e_{3}}e_{3}=-e_{1},%
\end{array}%
\text{ \ \ \ }%
\begin{array}{c}
\nabla _{e_{4}}e_{3}=e_{2},%
\end{array}
\\
&&%
\begin{array}{c}
\nabla _{e_{1}}e_{4}=e_{2},%
\end{array}%
\text{ \ \ \ \ }%
\begin{array}{c}
\nabla _{e_{2}}e_{4}=e_{1},%
\end{array}%
\text{ \ \ }%
\begin{array}{c}
\nabla _{e_{3}}e_{4}=e_{2},%
\end{array}%
\text{ \ \ \ \ \ \ }%
\begin{array}{c}
\nabla _{e_{4}}e_{4}=e_{1}.%
\end{array}%
\end{eqnarray*}%
The only non-zero components of the Riemann curvature tensor $R,$ with
respect to $\left\{ e_{1},e_{2},e_{3},e_{4}\right\} ,$ are 
\begin{equation*}
R_{1212}=R_{1214}=-R_{1223}=-R_{1234}=-R_{1434}=R_{2334}=-R_{3434}=-2.
\end{equation*}%
Consequently, the non-zero components of the Ricci tensor are given by%
\begin{eqnarray*}
&&%
\begin{array}{c}
\varrho \left( e_{1},e_{1}\right) =\varrho \left( e_{2},e_{2}\right)
=\varrho \left( e_{3},e_{3}\right) =\varrho \left( e_{4},e_{4}\right) =-2,%
\end{array}
\\
&&%
\begin{array}{c}
\varrho \left( e_{1},e_{3}\right) =\varrho \left( e_{2},e_{4}\right) =-4.%
\end{array}%
\end{eqnarray*}
\end{lm}

Next, let $\mathrm{D}\in \mathrm{Der}\left( \mathfrak{g}\right) $ where $%
\mathfrak{g}$ is the Lie algebra used in 
\eqref{LieB}
and put 
\begin{equation*}
\mathrm{D}U_{l}=\lambda _{l}^{1}U_{1}+\lambda _{l}^{2}U_{2}+\lambda
_{l}^{3}U_{3}+\lambda _{l}^{4}U_{4}+\lambda _{l}^{5}U_{5}\text{ \ for all }%
l=1,..,5.
\end{equation*}

Using 
\eqref{LieB}%
, we prove that 
\eqref{Der}
is satisfied if and only if%
\begin{eqnarray}
&&%
\begin{array}{c}
\lambda _{3}^{4}+2\left( \lambda _{1}^{2}-\lambda _{1}^{5}\right) =0,%
\end{array}
\label{sysDB} \\
&&%
\begin{array}{c}
\lambda _{1}^{5}+\varepsilon \left( \lambda _{3}^{4}-\lambda _{1}^{2}\right)
=0,%
\end{array}
\notag \\
&&%
\begin{array}{c}
\lambda _{1}^{2}-\lambda _{4}^{3}+2\lambda _{1}^{5}-\varepsilon \lambda
_{5}^{1}-\lambda _{2}^{1}=0,%
\end{array}
\notag \\
&&%
\begin{array}{c}
\lambda _{1}^{1}+\lambda _{4}^{4}-\varepsilon \lambda _{5}^{2}-\lambda
_{2}^{2}=0,%
\end{array}
\notag \\
&&%
\begin{array}{c}
\lambda _{2}^{5}-\varepsilon \left( \lambda _{1}^{1}+\lambda
_{4}^{4}-\lambda _{5}^{5}\right) =0,%
\end{array}
\notag \\
&&%
\begin{array}{c}
\lambda _{3}^{4}-2\lambda _{2}^{1}+\varepsilon \lambda _{5}^{1}=0,%
\end{array}
\notag \\
&&%
\begin{array}{c}
\lambda _{2}^{5}+\varepsilon \left( \lambda _{2}^{2}-\lambda _{5}^{5}\right)
=0,%
\end{array}
\notag \\
&&%
\begin{array}{c}
\lambda _{1}^{1}-\lambda _{2}^{2}-\lambda _{4}^{4}-2\lambda _{2}^{5}=0,%
\end{array}
\notag \\
&&%
\begin{array}{c}
\lambda _{1}^{2}-\lambda _{2}^{1}-\lambda _{4}^{3}=0,%
\end{array}
\notag \\
&&%
\begin{array}{c}
\lambda _{1}^{5}+\varepsilon \left( \lambda _{4}^{3}-\lambda _{2}^{1}\right)
=0,%
\end{array}
\notag
\end{eqnarray}

\begin{eqnarray*}
&&%
\begin{array}{c}
\lambda _{3}^{2}+\lambda _{4}^{1}+2\lambda _{3}^{5}=0,%
\end{array}
\notag \\
&&%
\begin{array}{c}
\lambda _{3}^{1}-\lambda _{4}^{2}+2\lambda _{4}^{5}=0,%
\end{array}
\notag \\
&&%
\begin{array}{c}
\lambda _{5}^{1}-2\left( \lambda _{2}^{1}+\lambda _{3}^{4}\right) =0,%
\end{array}
\notag \\
&&%
\begin{array}{c}
\lambda _{5}^{2}+2\left( \lambda _{2}^{2}-\lambda _{5}^{5}\right) =0,%
\end{array}
\notag \\
&&%
\begin{array}{c}
\lambda _{5}^{2}+2\left( -\lambda _{1}^{1}+\lambda _{4}^{4}+\lambda
_{5}^{5}\right) =0,%
\end{array}
\notag \\
&&%
\begin{array}{c}
\lambda _{5}^{1}-2\left( \lambda _{1}^{2}+\lambda _{4}^{3}\right) =0,%
\end{array}
\notag \\
&&%
\begin{array}{c}
\lambda _{4}^{2}=-\lambda _{3}^{1},\text{ \ }\lambda _{5}^{1}=2\varepsilon
\lambda _{1}^{5},\text{ \ }\lambda _{5}^{2}=2\varepsilon \lambda _{2}^{5},%
\text{ \ }%
\end{array}
\notag \\
&&%
\begin{array}{c}
\lambda _{1}^{3}=\lambda _{1}^{4}=\lambda _{2}^{3}=\lambda _{2}^{4}=\lambda
_{3}^{3}=\lambda _{5}^{3}=\lambda _{5}^{4}=0.%
\end{array}
\notag
\end{eqnarray*}%
So, we need to consider two cases:

\begin{itemize}
\item If $\varepsilon =1.$ In this case, we prove that all solutions of 
\eqref{sysDB}
are given by%
\begin{eqnarray*}
&&%
\begin{array}{c}
\lambda _{1}^{5}=\lambda _{2}^{1}=\lambda _{1}^{2},%
\end{array}%
\begin{array}{c}
\lambda _{4}^{1}=-\lambda _{3}^{2}-2\lambda _{3}^{5},%
\end{array}%
\begin{array}{c}
\lambda _{4}^{2}=\lambda _{4}^{5}=-\lambda _{3}^{1},%
\end{array}
\\
&&%
\begin{array}{c}
\lambda _{2}^{2}=\lambda _{1}^{1}-2\lambda _{2}^{5},%
\end{array}%
\begin{array}{c}
\lambda _{5}^{1}=2\lambda _{1}^{2},%
\end{array}%
\begin{array}{c}
\lambda _{5}^{2}=2\lambda _{2}^{5},%
\end{array}%
\begin{array}{c}
\lambda _{5}^{5}=\lambda _{1}^{1}-\lambda _{2}^{5},%
\end{array}
\\
&&%
\begin{array}{c}
\lambda _{3}^{4}=\lambda _{4}^{3}=\lambda _{4}^{4}=0.%
\end{array}%
\end{eqnarray*}

\item If $\varepsilon =-1,$ all solutions of 
\eqref{sysDB}
are given by%
\begin{eqnarray*}
&&%
\begin{array}{c}
\lambda _{4}^{1}=-\lambda _{3}^{2}-2\lambda _{3}^{5},%
\end{array}%
\begin{array}{c}
\lambda _{4}^{2}=\lambda _{4}^{5}=-\lambda _{3}^{1},%
\end{array}%
\begin{array}{c}
\lambda _{5}^{5}=\lambda _{2}^{2}=\lambda _{1}^{1},%
\end{array}
\\
&&%
\begin{array}{c}
\lambda _{1}^{2}=\lambda _{1}^{5}=\lambda _{2}^{1}=\lambda _{2}^{5}=\lambda
_{3}^{4}=\lambda _{4}^{3}=\lambda _{4}^{4}=\lambda _{5}^{1}=\lambda
_{5}^{2}=0.%
\end{array}%
\end{eqnarray*}
\end{itemize}

Therefore, we proved the following.

\begin{lemma}
Let $\mathfrak{g}=\mathfrak{h\oplus }m$ be the Lie algebra used in 
\eqref{LieB}%
. Then $\mathrm{D}\in \mathrm{Der}\left( \mathfrak{g}\right) $ if and only if

\begin{itemize}
\item $\varepsilon =1$%
\begin{equation*}
\mathrm{D}=\left( 
\begin{array}{ccccc}
\lambda _{1}^{1} & \lambda _{1}^{2} & \lambda _{3}^{1} & -\lambda
_{3}^{2}-2\lambda _{3}^{5} & 2\lambda _{1}^{2} \\ 
\lambda _{1}^{2} & \lambda _{1}^{1}-2\lambda _{2}^{5} & \lambda _{3}^{2} & 
-\lambda _{3}^{1} & 2\lambda _{2}^{5} \\ 
0 & 0 & 0 & 0 & 0 \\ 
0 & 0 & 0 & 0 & 0 \\ 
\lambda _{1}^{2} & \lambda _{2}^{5} & \lambda _{3}^{5} & -\lambda _{3}^{1} & 
\lambda _{1}^{1}-\lambda _{2}^{5}%
\end{array}%
\right)
\end{equation*}

\item $\varepsilon =-1$%
\begin{equation*}
\mathrm{D}=\left( 
\begin{array}{ccccc}
\lambda _{1}^{1} & 0 & \lambda _{3}^{1} & -\lambda _{3}^{2}-2\lambda _{3}^{5}
& 0 \\ 
0 & \lambda _{1}^{1} & \lambda _{3}^{2} & -\lambda _{3}^{1} & 0 \\ 
0 & 0 & 0 & 0 & 0 \\ 
0 & 0 & 0 & 0 & 0 \\ 
0 & \lambda _{2}^{5} & \lambda _{3}^{5} & -\lambda _{3}^{1} & \lambda
_{1}^{1}%
\end{array}%
\right) .
\end{equation*}
\end{itemize}
\end{lemma}

Using the above lemma, we now prove the following.

\begin{theorem}
Let $\left( M=G/H,g\right) $ be a four-dimensional generalized symmetric
space of type B. Then, $M$ is not an algebraic Ricci solitons.
\end{theorem}

\textbf{Proof. }Using Lemma \ref{B} we obtain that the Ricci operator of $%
\left( M=G/H,g\right) $ is given, with respect to the basis $\left\{
U_{1},U_{2},U_{3},U_{4},U_{5}\right\} ,$ by%
\begin{equation*}
\mathrm{Ric}=\left( 
\begin{array}{cccc}
-4\lambda & 0 & 4\lambda ^{2}+3 & 0 \\ 
0 & -4\lambda & 0 & 4\lambda ^{2}+3 \\ 
-4 & 0 & 4\lambda & 0 \\ 
0 & -4 & 0 & 4\lambda%
\end{array}%
\right) .
\end{equation*}%
Hence it follows, from the above lemma, that the algebraic Ricci soliton
condition 
\eqref{AS}
on $M$ does not occur. $\square $

\section{\protect\LARGE Algebraic Ricci soliton of spaces of type C }

\setcounter{equation}{0}

Let $\left( M=G,g\right) $ be a four-dimensional generalized symmetric space
of type C. Without loss of generality, we assume that the signature is $%
\left( 3,1\right) .$ The Lie algebra $\mathfrak{g}$ admits a basis $\left\{
U_{1},U_{2},U_{3},U_{4}\right\} $, such that the Lie bracket $\left[ \text{ }%
,\right] $ and the scalar product $\left\langle \text{ },\right\rangle $ on $%
\mathfrak{g}$ are given, respectively, by%
\begin{equation}
\begin{array}{c|cccc}
\left[ \text{ },\right] & U_{1} & U_{2} & U_{3} & U_{4} \\ \hline
U_{1} & 0 & 0 & 0 & -U_{1} \\ 
U_{2} & 0 & 0 & 0 & U_{2} \\ 
U_{3} & 0 & 0 & 0 & 0 \\ 
U_{4} & U_{1} & -U_{2} & 0 & 0%
\end{array}
\label{LieC}
\end{equation}%
and 
\begin{equation*}
\begin{array}{c|cccc}
\left\langle \text{ },\right\rangle & U_{1} & U_{2} & U_{3} & U_{4} \\ \hline
U_{1} & 1 & 0 & 0 & 0 \\ 
U_{2} & 0 & 1 & 0 & 0 \\ 
U_{3} & 0 & 0 & 0 & 1/2 \\ 
U_{4} & 0 & 0 & 1/2 & 0%
\end{array}%
\end{equation*}

The following result was proven in \cite{CD}.

\begin{lm}
\label{C} Let $M$ be a four-dimensional generalized symmetric space of type
C, with signature $\left( 3,1\right) .$ Then, there exist a
pseudo-orthonormal frame field 
\begin{equation*}
e_{1}=U_{1},\text{ \ \ }e_{2}=U_{2},\text{ \ \ }e_{3}=U_{3}+U_{4},\text{ \ \ 
}e_{4}=U_{3}-U_{4},
\end{equation*}%
on $M$, with $\left\langle e_{1},e_{1}\right\rangle =\left\langle
e_{2},e_{2}\right\rangle =\left\langle e_{3},e_{3}\right\rangle
=-\left\langle e_{4},e_{4}\right\rangle =1.$ The non-vanishing components of
the Levi-Civita connection $\nabla $ of $M$ are given by%
\begin{equation*}
\begin{array}{c}
\nabla _{e_{1}}e_{1}=-\nabla _{e_{2}}e_{2}=e_{3}+e_{4},%
\end{array}%
\text{ \ \ }%
\begin{array}{c}
\nabla _{e_{1}}e_{4}=-\nabla _{e_{1}}e_{3}=e_{1},%
\end{array}%
\text{ \ \ }%
\begin{array}{c}
\nabla _{e_{2}}e_{3}=-\nabla _{e_{2}}e_{4}=e_{2}.%
\end{array}%
\text{ }
\end{equation*}%
The non-zero components of the Riemann curvature tensor $R,$ with respect to 
$\left\{ e_{1},e_{2},e_{3},e_{4}\right\} ,$ are 
\begin{equation*}
R_{1313}=-R_{1314}=R_{1414}=R_{2323}=-R_{2324}=R_{2424}=-1.
\end{equation*}%
Consequently, the only non-zero components of the Ricci tensor are given by%
\begin{equation*}
\varrho \left( e_{3},e_{3}\right) =\varrho \left( e_{4},e_{4}\right)
=-\varrho \left( e_{3},e_{4}\right) =-2.
\end{equation*}
\end{lm}

Next, put $\mathrm{D}U_{l}=\lambda _{l}^{1}U_{1}+\lambda
_{l}^{2}U_{2}+\lambda _{l}^{3}U_{3}+\lambda _{l}^{4}U_{4}$ for all $%
l=1,..,4, $ where $\{U_{1},U_{2},U_{3},U_{4}\}$ is the basis used in 
\eqref{LieC}%
. Standard computations proves that $\mathrm{D}\in \mathrm{Der}\left( 
\mathfrak{g}\right) $ if and only if 
\begin{equation*}
\begin{array}{c}
\lambda _{1}^{2}=\lambda _{1}^{3}=\lambda _{1}^{4}=\lambda _{2}^{1}=\lambda
_{2}^{3}=\lambda _{2}^{4}=\lambda _{3}^{1}=\lambda _{3}^{2}=\lambda
_{3}^{4}=\lambda _{4}^{4}=0.%
\end{array}%
\end{equation*}%
So, we deduce the following.

\begin{lemma}
Let $\mathfrak{g}=\mathfrak{h\oplus }m$ be the Lie algebra used in 
\eqref{LieC}%
. Then $\mathrm{D}\in \mathrm{Der}\left( \mathfrak{g}\right) $ if and only if%
\begin{equation*}
\mathrm{D}=\left( 
\begin{array}{cccc}
\lambda _{1}^{1} & 0 & 0 & \lambda _{4}^{1} \\ 
0 & \lambda _{2}^{2} & 0 & \lambda _{4}^{2} \\ 
0 & 0 & \lambda _{3}^{3} & \lambda _{4}^{3} \\ 
0 & 0 & 0 & 0%
\end{array}%
\right) .
\end{equation*}
\end{lemma}

We can now prove the following.

\begin{theorem}
Let $\left( M=G/H,g\right) $ be a four-dimensional generalized symmetric
space of type C. Then, $M$ is an algebraic Ricci solitons. In particular%
\begin{equation*}
\mathrm{D=Ric=}\left( 
\begin{array}{cccc}
0 & 0 & 0 & 0 \\ 
0 & 0 & 0 & 0 \\ 
0 & 0 & 0 & -4 \\ 
0 & 0 & 0 & 0%
\end{array}%
\right) \text{ and }c=0.
\end{equation*}
\end{theorem}

\textbf{Proof. }Using Lemma \ref{C} we write down the Ricci operator of $%
\left( M=G/H,g\right) $, with respect to the basis $\left\{
U_{1},U_{2},U_{3},U_{4}\right\} ,$ getting%
\begin{equation*}
\mathrm{Ric}=\left( 
\begin{array}{cccc}
0 & 0 & 0 & 0 \\ 
0 & 0 & 0 & 0 \\ 
0 & 0 & 0 & -4 \\ 
0 & 0 & 0 & 0%
\end{array}%
\right) .
\end{equation*}%
Thus, using the above lemma, we obtain that the algebraic Ricci soliton
condition 
\eqref{AS}
on $M$ is satisfied if and only if 
\begin{equation*}
\lambda _{1}^{1}=\lambda _{2}^{2}=\lambda _{3}^{3}=\lambda _{4}^{1}=\lambda
_{4}^{2}=c=0\text{ and }\lambda _{4}^{3}=-4.\square
\end{equation*}

\section{\protect\LARGE Algebraic Ricci soliton of spaces of type D }

\setcounter{equation}{0}

Let $\left( M=G/H,g\right) $ be a four-dimensional generalized symmetric
space of type D, with signature $\left( 2,2\right) $. The Lie algebra $%
\mathfrak{g}=\mathfrak{h\oplus }m$ of the Lie group $G$ admits a basis $%
\left\{ U_{1},U_{2},U_{3},U_{4},U_{5}\right\} $, with $\left\{
U_{1},U_{2},U_{3},U_{4}\right\} $ and $\left\{ U_{5}\right\} $ are
respectively a basis of $m$ and of $\mathfrak{h}$, such that%
\begin{equation}
\begin{array}{c|ccccc}
\left[ \text{ },\right] & U_{1} & U_{2} & U_{3} & U_{4} & U_{5} \\ \hline
U_{1} & 0 & 0 & 0 & -U_{2} & U_{1} \\ 
U_{2} & 0 & 0 & -U_{1} & 0 & -U_{2} \\ 
U_{3} & 0 & U_{1} & 0 & -U_{5} & 2U_{3} \\ 
U_{4} & U_{2} & 0 & U_{5} & 0 & -2U_{4} \\ 
U_{5} & -U_{1} & U_{2} & -2U_{3} & 2U_{4} & 0%
\end{array}
\label{LieD}
\end{equation}%
and 
\begin{equation*}
\begin{array}{c|cccc}
\left\langle \text{ },\right\rangle & U_{1} & U_{2} & U_{3} & U_{4} \\ \hline
U_{1} & 0 & 1 & 0 & 0 \\ 
U_{2} & 1 & 0 & 0 & 0 \\ 
U_{3} & 0 & 0 & 0 & \lambda \\ 
U_{4} & 0 & 0 & \lambda & 0%
\end{array}%
\end{equation*}

with $\lambda \neq 0$ is a real constant.

The following result was proven in \cite{CD}.

\begin{lm}
\label{D} Let $M$ be a four-dimensional generalized symmetric space of type
D, with signature $\left( 2,2\right) .$ Then, there exist a
pseudo-orthonormal frame field 
\begin{eqnarray*}
e_{1} &=&\frac{1}{\sqrt{2}}\left( U_{1}+U_{2}\right) ,\text{ \ \ \ \ \ }%
e_{2}=\frac{1}{\sqrt{2\left\vert \lambda \right\vert }}\left(
U_{3}+\varepsilon U_{4}\right) , \\
e_{3} &=&\frac{1}{\sqrt{2}}\left( U_{1}-U_{2}\right) ,\text{ \ \ \ \ \ }%
e_{4}=\frac{1}{\sqrt{2\left\vert \lambda \right\vert }}\left(
U_{3}-\varepsilon U_{4}\right) ,
\end{eqnarray*}%
on $M$, with $\varepsilon =\pm 1$ and $\left\langle e_{1},e_{1}\right\rangle
=\left\langle e_{2},e_{2}\right\rangle =-\left\langle
e_{3},e_{3}\right\rangle =-\left\langle e_{4},e_{4}\right\rangle =1.$ The
non-vanishing components of the Levi-Civita connection $\nabla $ of $M$ are
given by%
\begin{eqnarray*}
&&%
\begin{array}{c}
\nabla _{e_{1}}e_{1}=\frac{1}{2\sqrt{2\left\vert \lambda \right\vert }}%
\left( \left( \varepsilon +1\right) e_{2}+\left( \varepsilon -1\right)
e_{4}\right) ,%
\end{array}
\\
&&%
\begin{array}{c}
\nabla _{e_{1}}e_{2}=-\frac{1}{2\sqrt{2\left\vert \lambda \right\vert }}%
\left( \left( \varepsilon +1\right) e_{1}-\left( \varepsilon -1\right)
e_{3}\right) ,%
\end{array}
\\
&&%
\begin{array}{c}
\nabla _{e_{1}}e_{3}=\frac{1}{2\sqrt{2\left\vert \lambda \right\vert }}%
\left( \left( \varepsilon -1\right) e_{2}+\left( \varepsilon +1\right)
e_{4}\right) ,%
\end{array}
\\
&&%
\begin{array}{c}
\nabla _{e_{1}}e_{4}=\frac{1}{2\sqrt{2\left\vert \lambda \right\vert }}%
\left( \left( \varepsilon -1\right) e_{1}-\left( \varepsilon +1\right)
e_{3}\right) ,%
\end{array}
\\
&&%
\begin{array}{c}
\nabla _{e_{3}}e_{1}=\frac{1}{2\sqrt{2\left\vert \lambda \right\vert }}%
\left( \left( \varepsilon -1\right) e_{2}+\left( \varepsilon +1\right)
e_{4}\right) ,%
\end{array}
\\
&&%
\begin{array}{c}
\nabla _{e_{3}}e_{2}=\frac{1}{2\sqrt{2\left\vert \lambda \right\vert }}%
\left( \left( 1-\varepsilon \right) e_{1}+\left( \varepsilon +1\right)
e_{3}\right) ,%
\end{array}
\\
&&%
\begin{array}{c}
\nabla _{e_{3}}e_{3}=\frac{1}{2\sqrt{2\left\vert \lambda \right\vert }}%
\left( \left( \varepsilon +1\right) e_{2}+\left( \varepsilon -1\right)
e_{4}\right) ,%
\end{array}
\\
&&%
\begin{array}{c}
\nabla _{e_{3}}e_{4}=\frac{1}{2\sqrt{2\left\vert \lambda \right\vert }}%
\left( \left( \varepsilon +1\right) e_{1}-\left( \varepsilon -1\right)
e_{3}\right) .%
\end{array}%
\end{eqnarray*}%
The non-zero components of the Riemann curvature tensor, with respect to $%
\left\{ e_{1},e_{2},e_{3},e_{4}\right\} ,$ are 
\begin{eqnarray*}
R_{1212} &=&-R_{1234}=-R_{1414}=-R_{1423}=-R_{2323}=R_{3434}=-\frac{1}{%
2\lambda }, \\
R_{1313} &=&-R_{1324}=-\frac{1}{\lambda }.
\end{eqnarray*}%
Consequently, the only non-zero components of the Ricci tensor are given by%
\begin{equation*}
\varrho \left( e_{2},e_{2}\right) =-\varrho \left( e_{4},e_{4}\right) =-%
\frac{1}{\lambda }.
\end{equation*}
\end{lm}

Now, let $\mathrm{D}\in \mathrm{Der}\left( \mathfrak{g}\right) $ where $%
\mathfrak{g}$ is the Lie algebra used in 
\eqref{LieD}%
. Put 
\begin{equation*}
\mathrm{D}U_{l}=\lambda _{l}^{1}U_{1}+\lambda _{l}^{2}U_{2}+\lambda
_{l}^{3}U_{3}+\lambda _{l}^{4}U_{4}+\lambda _{l}^{5}U_{5}\text{ \ for all }%
l=1,..,5.
\end{equation*}

Starting from 
\eqref{LieD}%
, we can write down 
\eqref{Der}
and we get%
\begin{eqnarray*}
&&%
\begin{array}{c}
\lambda _{3}^{3}=\lambda _{1}^{1}-\lambda _{2}^{2}\text{ \ },\lambda
_{3}^{5}=\lambda _{1}^{2},\text{\ }\lambda _{4}^{4}=-\lambda _{3}^{3},\text{
\ }\lambda _{4}^{5}=-\lambda _{2}^{1},%
\end{array}
\\
&&%
\begin{array}{c}
\text{ }\lambda _{5}^{1}=-\lambda _{4}^{2},\text{ }\lambda _{5}^{2}=\lambda
_{3}^{1},\text{ }\lambda _{5}^{3}=2\lambda _{2}^{1},\text{ }\lambda
_{5}^{4}=-2\lambda _{1}^{2},%
\end{array}
\\
&&%
\begin{array}{c}
\lambda _{1}^{3}=\lambda _{1}^{4}=\lambda _{1}^{5}=\lambda _{2}^{3}=\lambda
_{2}^{4}=\lambda _{2}^{5}=\lambda _{3}^{2}=\lambda _{3}^{4}=\lambda
_{4}^{1}=\lambda _{4}^{3}=\lambda _{5}^{5}=0.%
\end{array}%
\end{eqnarray*}%
We deduce the following.

\begin{lemma}
Let $\mathfrak{g}=\mathfrak{h\oplus }m$ be the Lie algebra used in 
\eqref{LieD}%
. Then $\mathrm{D}\in \mathrm{Der}\left( \mathfrak{g}\right) $ if and only if%
\begin{equation*}
\mathrm{D}=\left( 
\begin{array}{ccccc}
\lambda _{1}^{1} & \lambda _{2}^{1} & \lambda _{3}^{1} & 0 & -\lambda
_{4}^{2} \\ 
\lambda _{1}^{2} & \lambda _{2}^{2} & 0 & \lambda _{4}^{2} & \lambda _{3}^{1}
\\ 
0 & 0 & \lambda _{1}^{1}-\lambda _{2}^{2} & 0 & 2\lambda _{2}^{1} \\ 
0 & 0 & 0 & \lambda _{2}^{2} -\lambda _1^1 & -2\lambda _{1}^{2} \\ 
0 & 0 & \lambda _{1}^{2} & -\lambda _{2}^{1} & 0%
\end{array}%
\right) .
\end{equation*}
\end{lemma}

Using the above lemma, we prove the following.

\begin{theorem}
Let $\left( M=G/H,g\right) $ be a four-dimensional generalized symmetric
space of type D. Then, $M$ is an algebraic Ricci solitons. In particular%
\begin{equation*}
pr\circ \mathrm{D=}\left( 
\begin{array}{cccc}
\frac{1}{\lambda } & 0 & 0 & 0 \\ 
0 & \frac{1}{\lambda } & 0 & 0 \\ 
0 & 0 & 0 & 0 \\ 
0 & 0 & 0 & 0%
\end{array}%
\right) \text{ and }c=-\frac{1}{\lambda }.
\end{equation*}
\end{theorem}

\textbf{Proof. }Using Lemma \ref{D} we write down the Ricci operator of $%
\left( M=G/H,g\right) $, with respect to the basis $\left\{
U_{1},U_{2},U_{3},U_{4}\right\} ,$ getting%
\begin{equation*}
\mathrm{Ric}=\left( 
\begin{array}{cccc}
0 & 0 & 0 & 0 \\ 
0 & 0 & 0 & 0 \\ 
0 & 0 & -\frac{1}{\lambda } & 0 \\ 
0 & 0 & 0 & -\frac{1}{\lambda }%
\end{array}%
\right) .
\end{equation*}%
Using the above lemma, we obtain that the algebraic Ricci soliton condition 
\eqref{AS}
on $M$ is satisfied if and only if%
\begin{equation*}
\lambda _{1}^{1}=\lambda _{2}^{2}=-c=\frac{1}{\lambda }\text{ and }\lambda
_{1}^{2}=\lambda _{2}^{1}=\lambda _{3}^{1}=\lambda _{3}^{3}=\lambda
_{4}^{2}=0.\square
\end{equation*}


\begin{thebibliography}{99}
\bibitem{BO} W. Batat and K. Onda, Algebraic Ricci Solitons of
three-dimensional Lorentzian Lie groups, arXiv : 1112.2455.

\bibitem{C07} G. Calvaruso, Homogeneous structures on three-dimensional
Lorentzian manifolds, \textit{J. Geom. Phys}. \textbf{57} (2007), no. 4,
1279--1291.

\bibitem{CD} G. Calvaruso and B. De Leo, Curvature properties of
four-dimensional generalized symmetric spaces, \textit{J. Geom}. \textbf{90}
(2008), 30--46.

\bibitem{CK04} B. Chow and D. Knopf, The Ricci flow: an introduction,
Mathematical Surveys and Monographs, vol. 110, American Mathematical
Society, Providence, RI, 2004.

\bibitem{CK} J. Cerny and O. Kowalski, Classification of generalized
symmetric pseudo-Riemannian spaces of dimension n 4, \textit{Tensor N. S}. 
\textbf{38} (1982), 256--267.

\bibitem{DM} B. De Leo and R. A. Marinosci, Homogeneous geodesics of
four-dimensional generalized symmetric pseudo-Riemannian spaces, \textit{%
Publ. Math. Debrecen} \textbf{73} (2008), 341--360.

\bibitem{FR} M.E. Fels and A.G. Renner, Non-reductive homogeneous
pseudo-Riemannian manifolds of dimension four, \textit{Canad. J. Math}. (2) 
\textbf{58} (2006), 282--311.

\bibitem{Jac} M. Jablonski, Homogeneous Ricci solitons, arXiv:1109.6556.

\bibitem{L01} J. Lauret, Ricci soliton homogeneous nilmanifolds, \textit{%
Math. Ann}. \textbf{319} (2001), no. 4, 715--733.

\bibitem{O11} K. Onda, Examples of sol-solitons in pseudo-Riemannian case,
arXiv : 1112.0424v2.

\bibitem{O} B. O'Neill, Semi-Riemannian Geometry. New York: Academic
Press,1983.
\end{thebibliography}
\end{document}